\documentclass[12pt]{article}

\usepackage{amsmath, epsfig, cite}
\usepackage{amssymb}
\usepackage{amsfonts}
\usepackage{latexsym}
\usepackage{graphicx}
\usepackage{extarrows}

\makeatletter

\newcommand{\Rmnum}[1]{\expandafter\@slowromancap\romannumeral #1@}
\makeatother

\newtheorem{thm}{Theorem}[section]

\newtheorem{lem}[thm]{Lemma}

\newcommand{\qed}{{\hfill\rule{4pt}{7pt}}}
\def\pf{\noindent {\it Proof.} }
\def\R{{\widetilde{R}}}

\numberwithin{equation}{section}

\makeatletter \@addtoreset{equation}{section} \makeatother

\begin{document}
\rule{0cm}{1cm}
\begin{center}
{\Large\bf Combinatorial Proof of the Inversion Formula on

the Kazhdan-Lusztig $R$-Polynomials}
\end{center}
 \vskip 6mm
 \begin{center}
{\small William Y.C. Chen$^1$, Neil J.Y. Fan$^2$, Alan  J.X. Guo$^3$, Peter L. Guo$^4$\\
Harry H.Y. Huang$^5$, Michael X.X. Zhong$^6$}

\vskip 4mm
$^{1,3,4,5,6}$Center for Combinatorics, LPMC-TJKLC\\
Nankai University,
Tianjin 300071,
P.R. China \\[3mm]

$^2$Department of Mathematics\\
Sichuan University, Chengdu, Sichuan 610064, P.R. China

\vskip 4mm

$^1$chen@nankai.edu.cn,  $^2$fan@scu.edu.cn, $^3$aalen@mail.nankai.edu.cn\\
$^4$lguo@nankai.edu.cn,   $^5$haoyanghuang@mail.nankai.edu.cn\\
$^6$michaelzhong@mail.nankai.edu.cn
  \end{center}

\begin{abstract}
Let $W$ be a Coxeter group, and for  $u,v\in W$,
 let $R_{u,v}(q)$ be the Kazhdan-Lusztig $R$-polynomial indexed by $u$ and $v$.
In this paper, we present a combinatorial proof
 of the inversion formula on $R$-polynomials due to Kazhdan and Lusztig.
This problem was raised by Brenti.
Based on Dyer's combinatorial interpretation of the $R$-polynomials in terms of
 increasing Bruhat paths, we reformulate  the inversion formula
in terms of $V$-paths. By a $V$-path from $u$ to $v$  with bottom $w$
we mean a pair $(\Delta_1,\Delta_2)$ of Bruhat paths such that
$\Delta_1$ is a decreasing  path from $u$ to $w$ and
$\Delta_2$ is an increasing path from $w$
to $v$. We find a reflection principle on $V$-paths,
which leads to a combinatorial
proof of the inversion formula. Moreover, we give two applications of the  reflection principle.
First, we restrict this involution  to $V$-paths from $u$ to $v$
with maximal length. This
 provides a direct interpretation for the equi-distribution property
that any nontrivial  interval  $[u,v]$ has as many elements of even length as
elements of odd length. This property
was obtained by Verma in his derivation of
the M\"obius function of the Bruhat order.
Second, using the reflection principle for the symmetric group,
we obtain a refinement of the inversion formula by restricting the
summation to permutations ending with a given element.
\end{abstract}



\section{Introduction}
Let $W$ be a Coxeter group. For $u, v\in W$,  let $R_{u,v}(q)$ be
the Kazhdan-Lusztig $R$-polynomial indexed by $u$ and $v$.
The following inversion
formula was obtained by Kazhdan and Lusztig \cite{K-L}:
\begin{equation}\label{1}
\sum_{u\leq w \leq v}(-1)^{\ell(w)-\ell(u)} R_{u,w}(q)R_{w,v}(q)=\delta_{u,v},
\end{equation}
where $\leq $ is the Bruhat order and $\ell$ is the length function.
The aim  of this paper is
to present a combinatorial interpretation of this formula.
This problem was raised by Brenti \cite{Brenti}.

To give a combinatorial proof of the above inversion formula,
we start with  Dyer's combinatorial description
 of the $R$-polynomials in terms of
increasing Bruhat paths
\cite{Dyer}.
Then we reformulate the inversion formula
in terms of $V$-paths.
For $u\leq w\leq v$, by a $V$-path from $u$ to $v$  with bottom $w$
we mean a pair $(\Delta_1,\Delta_2)$ of Bruhat paths such that
$\Delta_1$ is a decreasing  path from $u$ to $w$ and
$\Delta_2$ is an increasing path from $w$
to $v$.
We construct  a reflection principle on $V$-paths. This leads to a
combinatorial proof of the inversion formula.

We give two applications of the reflection principle.
First, we restrict the reflection principle
to $V$-paths from $u$ to $v$
with maximal length. This induces an involution on the interval $[u,v]$ with $u<v$,
which leads to  a combinatorial proof of
  the equi-distribtution property that
 any nontrivial  interval  $[u,v]$ has as many elements of even length as
elements of odd length.
This property was first conjectured by Verma \cite{Verma2}, and
proved in \cite{Verma1} by induction. It has been used in
Verma's derivation of   the
M\"obius function of the Bruhat order \cite{Verma1}.
Other
 proofs of the M\"obius function formula for the Bruhat order can be found in \cite{BW,Deodhar,Stembridge}.
 Recently, Jones \cite{Jones} found a combinatorial proof
for the equi-distribution property by constructing an involution on the intervals
of a Coxeter group $W$. When $W$ is finite,
Jones \cite{Jones} showed that this involution agrees with the construction of
Rietsch and Williams \cite{Williams} in their study
 of discrete Morse theory and totally nonnegative flag varieties.
Our involution seems to be simpler as far as the equi-distribution property is
concerned.

As a second application,
we  find a refinement of  the inversion formula when  $W$ is the symmetric group
$S_n$. For a permutation $w\in S_n$, we write $w=w(1)w(2)\cdots w(n)$, where,
 for $1\leq i\leq n$,
 $w(i)$ denotes the element in
the $i$-th position.
Let $u$ and $v$ be two permutations in $S_n$
 such that $u<v$ in the Bruhat order. For  $1\leq k\leq n$,
let $[u,v]_k$ denote the set of permutations in the interval $[u,v]$
that end with $k$, that is,
\[[u,v]_k=\{w\in [u,v]\colon w(n)=k\}.\]
Applying  a variation of the reflection principle,
we show that the summation
\[\sum_{w\in [u,v]_k}(-1)^{\ell(w)-\ell(u)}R_{u,w}(q)R_{w,v}(q)\]
equals zero or a power of $q$ up to a sign.

\section{Preliminaries}

In this section, we give an overview of the terminology and background
 on Coxeter groups and $R$-polynomials that will be used in this paper.

Let $(W,S)$ be a Coxeter system, and
\[T=\{wsw^{-1}: s\in S, \ w\in W\}\]
its set of reflections.  We use $\ell(w)$ to denote the length of $w\in W$.
For two elements $u$ and $v$ in $W$,
$u\leq v$ in the Bruhat order if and only if there exists
a sequence of reflections $t_1, t_2,\ldots, t_r$  in $T$ such that
\begin{itemize}
\item[(i)] $v=u\,t_1\,t_2 \cdots  t_r$;
\item[(ii)] $\ell(u \,t_1\cdots t_i )>\ell( u \,t_1\cdots t_{i-1})$  for $1\leq i\leq r$.
\end{itemize}

The $R$-polynomials were  introduced by Kazhdan and Lusztig \cite{K-L} in order to describe  the
 inverses of the basis elements in  the Hecke algebra associated to
the Coxeter system $(W,S)$. They   obtained the following
recurrence relation. For $w\in W$,
let $D_R(w)$ be the right descent set of $w$, that is,
\[D_R(w)=\{s\in S\colon
\ell(ws)<\ell(w)\}.\]

\begin{thm}\label{pre1}
For any $u,v\in W$,
\begin{itemize}
\item[$\mathrm{(i)}$] $R_{u,v}(q)=0$, if $u\nleq v$;
\item[$\mathrm{(ii)}$] $R_{u,v}(q)=1$, if $u= v$;
\item[$\mathrm{(iii)}$] If $u<v$ and $s\in D_{R}(v)$, then
\[R_{u,v}(q)=\left\{
        \begin{array}{ll}
          R_{us,\,vs}(q), & \hbox{\rm{if} $s\in D_{R}(u)$;} \\[5pt]
          qR_{us,\,vs}(q)+(q-1)R_{u,\,vs}(q), & \hbox{\rm{if} $s \notin  D_{R}(u)$.}
        \end{array}
      \right.
\]
\end{itemize}
\end{thm}

Notice that the $R$-polynomials may contain negative coefficients.
A variant of the $R$-polynomials introduced by Dyer \cite{Dyer}, which
 have been called the $\widetilde{R}$-polynomials,
have nonnegative coefficients.
The following two theorems are due to \cite{Dyer}, see also
Bj\"orner and Brenti \cite{Bjorner}.

\begin{thm}\label{Rtilde}
For any $u,v\in W$,
\begin{itemize}
\item[$\mathrm{(i)}$] $\widetilde{R}_{u,v}(q)=0$, if $u\nleq v$;
\item[$\mathrm{(ii)}$] $\widetilde{R}_{u,v}(q)=1$, if $u= v$;
\item[$\mathrm{(iii)}$] If $u<v$ and $s\in D_{R}(v)$, then
\[\widetilde{R}_{u,v}(q)=\left\{
        \begin{array}{ll}
          \widetilde{R}_{us,\,vs}(q), & \hbox{\rm{if} $s\in D_{R}(u)$;} \\[5pt]
          \widetilde{R}_{us,\,vs}(q)+q\widetilde{R}_{u,\,vs}(q), & \hbox{\rm{if} $s \notin  D_{R}(u)$.}
        \end{array}
      \right.
\]
\end{itemize}
\end{thm}

\begin{thm}\label{avn}
For  $u\leq v\in W$, we have
\begin{equation}\label{f}
R_{u,v}(q)
=q^{\frac{\ell(v)-\ell(u)}{2}}\widetilde{R}_{u,v}(q^{\frac{1}{2}}-q^{-\frac{1}{2}}).
\end{equation}
\end{thm}

In this paper, we shall give a combinatorial interpretation of
the inversion formula for  $\R$-polynomials.
By Theorem \ref{avn}, the inversion formula \eqref{1} can
be restated as follows
\begin{equation}\label{3-r}
\sum_{u\leq w \leq v}(-1)^{\ell(w)-\ell(u)} \R_{u,w}(q)\R_{w,v}(q)=\delta_{u,v}.
\end{equation}

Dyer \cite{Dyer} has given a combinatorial interpretation of
the $\R$-polynomials in terms of increasing paths
in the Bruhat graph  of a Coxeter group.
The Bruhat graph $BG(W)$ of a Coxeter group $W$ is a
directed graph whose nodes are the elements of $W$
such that there is an arc from $u$ to $v$
if $v=u\,t$ for some reflection $t\in T$ and
$\ell(u)<\ell(v)$.

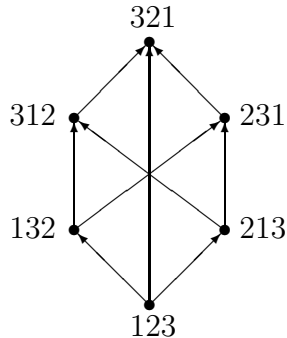
\begin{figure}[!htb]
\setlength{\unitlength}{2mm}
\begin{center}
  \begin{picture}(25,25)
    \multiput(10,3)(0,17.5){2}{\circle*{0.7}}
    \multiput(5,8)(0,7.5){2}{\circle*{0.7}}
    \multiput(15,8)(0,7.5){2}{\circle*{0.7}}
    \put(10,3){\vector(-1,1){4.75}}
    \put(10,3){\vector(1,1){4.75}}
    \multiput(5,8)(10,0){2}{\vector(0,1){7.25}}
    \put(5,15.5){\vector(1,1){4.75}}
    \put(15,15.5){\vector(-1,1){4.75}}
    \put(5,8){\vector(4,3){9.7}}
    \put(15,8){\vector(-4,3){9.7}}
    \put(10,3){\vector(0,1){17.25}}
    \put(8.7,1){$123$}
    \put(0.7,7.5){$132$}
    \put(16,7.5){$213$}
    \put(16,15){$231$}
    \put(0.7,15){$312$}
    \put(8.7,21.5){$321$}
  \end{picture}
\end{center}
\caption{The Bruhat graph of $S_3$}\label{fig}
\end{figure}

The Bruhat graph of  the symmetric group $S_n$ can be described as follows.
Let $(i,j)$ denote a  transposition in $S_n$,
where $1\leq i<j\leq n$.
The reflection set $T$ of $S_n$
consists of the
transpositions of $S_n$, that is,
\[T=\{(i,j)\colon 1\leq i<j\leq n\}.\]
For two permutations $u,v$ in $S_n$,
there is an arc from $u$ to $v$
if $v=u(i,j)$ and $u(i)<u(j)$, see \cite{Bjorner}.
The Bruhat graph of $S_3$ is illustrated in Figure \ref{fig}.

An increasing path in the Bruhat graph is defined based on
 the reflection ordering on
the positive roots  of a Coxeter group.
Let $\Phi$ be the root system  of $W$, and
 $\Phi^{+}$ be the positive root system.
A total ordering $\prec$ on $\Phi^{+}$
is called a reflection ordering on $\Phi^{+}$
 if for any $\alpha\prec\beta\in\Phi^+$ and two nonnegative  real numbers $\lambda,\mu$ such that $\lambda\alpha+\mu\beta\in\Phi^+$, we have
\[\alpha\prec\lambda\alpha+\mu\beta\prec\beta.\]

Here is an example of  a reflection ordering
for  the symmetric group $S_n$. Consider the following root system
\[\Phi=\{\pm(\varepsilon_i-\varepsilon_j)\colon 1\leq i< j \leq n\},\]
where $\varepsilon_i$ is the $i$-th coordinate vector,
and consider the following positive root system
\[\Phi^{+}=\{\varepsilon_i-\varepsilon_j\colon 1\leq i< j \leq n\}.\]
It is easy to verify that the following lexicographic ordering on $\Phi^{+}$ is a reflection ordering
\begin{equation}\label{2}
\varepsilon_1-\varepsilon_2\prec\varepsilon_1-\varepsilon_3\prec
\cdots<\varepsilon_1-\varepsilon_n\prec\varepsilon_2-\varepsilon_3
\prec\cdots\prec\varepsilon_{n-1}-\varepsilon_n.
\end{equation}

Since positive roots in $\Phi^+$ are in one-to-one
correspondence with reflections,
a reflection ordering induces  a total ordering on the set
$T$ of reflections.
For example, for the symmetric group $S_n$, a positive root $\varepsilon_i-\varepsilon_j$
corresponds to a transposition $(i,j)$. Thus,
the total ordering on transpositions
inheriting from the lexicographic ordering in \eqref{2}
is
\begin{equation}\label{lexico}
(1,2)\prec(1,3)\prec\cdots\prec(1,n)\prec(2,3)\prec\cdots\prec(n-1,n).
\end{equation}

We use $u\rightarrow v$ to denote an arc from $u$ to $v$ in the Bruhat graph.
Let $\Delta=u_0\rightarrow u_1\rightarrow\cdots \rightarrow u_r$ be a path
from $u$ to $v$, where $u_0=u$ and $u_r=v$.
We say $\Delta$ is
an increasing path with respect to a reflection ordering $\prec$ if
\[t_1\prec t_2\prec\cdots\prec t_r,\]
where, for $1\leq i\leq r$, $t_i$ is the reflection determined by $u_i=u_{i-1}t_i$.
For example, the following is an increasing path
with respect to the lexicographic ordering in the Bruhat graph of $S_4$
\[{2314}\xlongrightarrow{(12)} 3214\xlongrightarrow{(14)} 4213\xlongrightarrow{(24)} 4312.\]
Similarly, we say that $\Delta=u_0\rightarrow u_1\rightarrow\cdots \rightarrow u_r$ is
a decreasing path with respect to a reflection ordering $\prec$ if
\[t_1\succ t_2\succ\cdots\succ t_r,\]
where, for $1\leq i\leq r$, $t_i$ is the reflection
determined by $u_i=u_{i-1}t_i$.

Let $\ell(\Delta)$ denote the length of $\Delta$, namely, the number of
arcs in $\Delta$.
The following combinatorial interpretation of the $\widetilde{R}$-polynomials
was obtained by Dyer \cite{Dyer}.

\begin{thm}\label{DyerComInterpretation}
Let $W$ be a Coxeter group, and $T$ its set of reflections.
Then, for
 any fixed reflection ordering $\prec$ on $T$ and any $u,v\in W$ with $u\leq v$,  we have
\begin{equation}\label{c}
\widetilde{R}_{u,v}(q)=\sum_{\Delta} q^{\ell(\Delta)},
\end{equation}
where the sum ranges over increasing Bruhat paths from $u$ to $v$ with respect to $\prec$.
\end{thm}

Since the reverse of
a reflection ordering  is also a reflection ordering,
relation \eqref{c} can be restated as follows
\begin{equation}\label{d}
\widetilde{R}_{u,v}(q)=\sum_{\Delta'} q^{\ell(\Delta')},
\end{equation}
where the sum ranges  over decreasing Bruhat paths from $u$ to $v$
 with respect
to any fixed reflection ordering. Both \eqref{c} and \eqref{d} will be
used in our combinatorial interpretation of the inversion formula  \eqref{3-r}  in the next section.

\section{A reflection principle  on $V$-paths}

In this section, we give a restatement of the inversion formula
in terms of $V$-paths, and we construct a reflection principle  on $V$-paths
 leading to a combinatorial proof of the inversion formula.

By a $V$-path  from $u$ to $v$ with bottom $w$ under a reflection ordering $\prec$,
we mean a
pair $(\Delta_1, \Delta_2)$ of Bruhat paths such that
$\Delta_1$ is a decreasing path from $u$
 to $w$ and
$\Delta_2$ is an increasing path from $w$ to $v$.
The sign of a $V$-path $(\Delta_1,\Delta_2)$ is defined by
\[\mathrm{sgn}(\Delta_1,\Delta_2)=(-1)^{\ell(\Delta_1)}.\]
It is known that the length  of a  Bruhat path  from $u$ to $w$ has the same parity
as $\ell(w)-\ell(u)$, see Bj\"orner and Brenti\cite{Bjorner}. It follows that
\begin{equation}\label{sign}
\mathrm{sgn}(\Delta_1,\Delta_2)=(-1)^{\ell(w)-\ell(u)}.
\end{equation}

Using \eqref{c}, \eqref{d} and  \eqref{sign},
we are led to
  the following restatement of the inversion formula.

\begin{thm}\label{thm-4}
For a Coxeter group $W$ and  $u\leq v\in W$, we have
\begin{equation}\label{3}
\sum_{u\leq w \leq v}(-1)^{\ell(w)-\ell(u)} \R_{u,w}(q)\R_{w,v}(q)=
\sum_{(\Delta_1,\Delta_2)}\mathrm{sgn}(\Delta_1,\Delta_2) q^{\ell(\Delta_1)+\ell(\Delta_2)}=\delta_{u,v},
\end{equation}
where the sum ranges over $V$-paths from $u$ to $v$ under any fixed reflection ordering.
\end{thm}

We shall construct
a reflection principle  $I$ on $V$-paths from $u$ to $v$,
where $u<v$ in the Bruhat order, which turns out to be
a length preserving and sign reversing involution. So we obtain
a combinatorial interpretation of the inversion formula.

\noindent
\textbf{Reflection Principle on $V$-Paths:}
For $u<v$, let $(\Delta_1,\Delta_2)$ be a $V$-path from $u$ to $v$ with
bottom   $w$. Write
\[\Delta_1=u_0\rightarrow u_1\rightarrow\cdots \rightarrow u_i\ \ \text{and}\ \
\Delta_2=v_0\rightarrow v_1\rightarrow\cdots \rightarrow v_j,\]
where $u_0=u$, $u_i=v_0=w$ and $v_j=v$.
From $(\Delta_1,\Delta_2)$, we construct a $V$-path $(\Delta_1',\ \Delta_2')$
by the following procedure:
\begin{itemize}
\item[(1)] $u=w$, that is, $\Delta_1$ is an empty path. Set
\[\Delta_1'=u\rightarrow v_1\ \ \text{and}\ \ \Delta_2'=v_1\rightarrow v_2\rightarrow\cdots
\rightarrow v_j.\]

\item[(2)] $v=w$, that is, $\Delta_2$ is an empty path. Set
\[\Delta_1'=u_0\rightarrow u_1\rightarrow\cdots\rightarrow u_{i-1}\ \ \text{and}\ \ \Delta_2'=u_{i-1}\rightarrow  v.\]

\item[(3)] $u\neq w$ and $v\neq w$. In this case, we may assume that
$u_r=u_{r-1}\,t_r$ for $1\leq r\leq i$, and
 $v_s=v_{s-1}\,t_s'$ for $1\leq s\leq j$. If $t_i\prec t_1'$, set
\[\Delta_1'=u_0\rightarrow u_1\rightarrow \cdots\rightarrow u_{i-1}\ \ \text{and}\ \ \Delta_2'=u_{i-1}\rightarrow v_0\rightarrow v_1\rightarrow\cdots\rightarrow v_j.\]
If $t_i\succ t_1'$, set
\[\Delta_1'=u_0\rightarrow u_1\rightarrow \cdots\rightarrow u_i\rightarrow v_1\ \ \text{and}\ \ \Delta_2'= v_1\rightarrow \cdots\rightarrow v_j.\]
\end{itemize}
The $V$-path $(\Delta_1', \Delta_2')$ is defined to be
$I(\Delta_1, \Delta_2)$.

\begin{figure}[!htb]
\setlength{\unitlength}{2mm}
\begin{center}
\begin{picture}(100,21)


\put(0,0){
  \begin{picture}(35,21)
    \multiput(19,3)(3,4){5}{\circle*{0.7}}
    \multiput(19,3)(3,4){2}{\vector(3,4){2.8}}
    \multiput(25,11)(0.3,0.4){10}{\circle*{0.1}}
    \put(28,15){\vector(3,4){2.8}}
    \multiput(19,3)(-3,4){5}{\circle*{0.7}}
    \multiput(13,11)(3,-4){2}{\vector(3,-4){2.8}}
    \multiput(13,11)(-0.3,0.4){10}{\circle*{0.1}}
    \put(7,19){\vector(3,-4){2.8}}
    \put(13.75,1){$u_i=w=v_0$}
    \put(11.75,6.5){$u_{i-1}$}
    \put(23,6.5){$v_{1}$}
    \put(8.75,10.5){$u_{i-2}$}
    \put(26,10.5){$v_{2}$}
    \put(7.5,14.5){$u_{1}$}
    \put(29,14.5){$v_{j-1}$}
    \put(0.5,18.5){$u=u_{0}$}
    \put(32,18.5){$v_j=v$}
    \put(16,4.5){$t_{i}$}
    \put(18.25,4.5){$\succ$}
    \put(21,4.5){$t'_{1}$}
    \put(14.75,9){$t_{i-1}$}
    \put(8.75,17){$t_1$}
    \put(21.5,9){$t'_2$}
    \put(27.5,17){$t'_{j}$}
  \end{picture}
}


\put(35,10){
  \begin{picture}(10,5)
    \put(0,2){\vector(1,0){8}}
    \put(8,0){\vector(-1,0){8}}
  \end{picture}
}


\put(40,0){
  \begin{picture}(35,21)
    \multiput(19,3)(3,4){5}{\circle*{0.7}}
    \multiput(19,3)(3,4){2}{\vector(3,4){2.8}}
    \multiput(25,11)(0.3,0.4){10}{\circle*{0.1}}
    \put(28,15){\vector(3,4){2.8}}
    \multiput(19,3)(-3,4){5}{\circle*{0.7}}
    \multiput(13,11)(3,-4){2}{\vector(3,-4){2.8}}
    \multiput(13,11)(-0.3,0.4){10}{\circle*{0.1}}
    \put(7,19){\vector(3,-4){2.8}}
    \put(11.9,1){$u_{i+1}=w'=v_1$}
    \put(9,6.5){$v_0=u_{i}$}
    \put(23,6.5){$v_{2}$}
    \put(8.75,10.5){$u_{i-1}$}
    \put(26,10.5){$v_{3}$}
    \put(7.5,14.5){$u_{1}$}
    \put(29,14.5){$v_{j-1}$}
    \put(0.5,18.5){$u=u_{0}$}
    \put(32,18.5){$v_j=v$}
    \put(16,4){$t'_{1}$}
    \put(18.25,4){$\prec$}
    \put(21,4){$t'_{2}$}
    \put(14.75,9){$t_{i}$}
    \put(8.75,17){$t_1$}
    \put(21.5,9){$t'_3$}
    \put(27.5,17){$t'_{j}$}
  \end{picture}
}
\end{picture}
\end{center}
\caption{The reflection principle on $V$-paths.}\label{fig-2}
\end{figure}
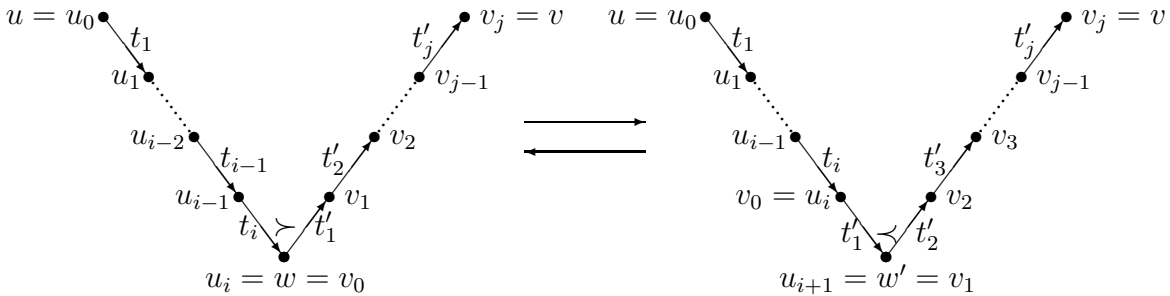

An illustration of the reflection principle is given in Figure \ref{fig-2}.
Let us consider an example.
Let $W$ be the symmetric group $S_4$. We choose the lexicographic order
in (\ref{lexico}) as the reflection ordering. Let $u=1234$ and $v=4312$.
The following is a
 $V$-path from
$u$ to $v$
\[1234\xlongrightarrow{(23)} 1324\xlongrightarrow{(13)} \underline{2314}\xlongrightarrow{(12)} 3214\xlongrightarrow{(14)} 4213\xlongrightarrow{(24)} 4312,\]
where we use  underline to indicate the bottom node.
Since $(1,3)\succ (1,2)$, we see that under the reflection principle $I$,
 this $V$-path is mapped to the $V$-path
\[1234\xlongrightarrow{(23)} 1324\xlongrightarrow{(13)} 2314\xlongrightarrow{(12)} \underline{3214}\xlongrightarrow{(14)} 4213\xlongrightarrow{(24)} 4312.\]

To conclude this section, we remark that the involution $I$
also yields a combinatorial interpretation of the following
result of Verma \cite{Verma1}.

\begin{thm}\label{verma-2}
Let $W$ be a Coxeter group and  $u<v\in W$. Then
the interval $[u,v]$ has the same number of elements of odd length
as elements of even length.
\end{thm}

{\noindent \it Combinatorial Proof.}
It is known that for $\sigma\leq \tau\in W$,
there exists a unique maximal increasing (or, decreasing)
Bruhat path  from $\sigma$ to $\tau$ \cite{Dyer}.
Thus, for any $w\in [u,v]$ there is a unique maximal $V$-path
from $u$ to $v$ whose bottom is $w$. So the maximal $V$-paths from $u$ to $v$ are in
one-to-one correspondence
with  elements in the interval $[u,v]$.

By the construction of the involution $I$, it can be seen that
the set  of all
maximal $V$-paths from $u$ to $v$ is invariant under $I$.
Thus, the restriction of   $I$
to the maximal $V$-paths from $u$ to $v$ induces an involution on
the elements in the interval
$[u,v]$. Moreover, the induced involution on the interval $[u,v]$
changes parity of the length of each element.
This completes the proof. \qed

\section{A refinement of the inversion formula for $S_n$}

In this section, we use the reflection principle to give  a refinement
of the inversion formula for the symmetric group $S_n$.
For $u<v\in S_n$ and $1\leq k \leq n$, let  $[u,v]_k$
denote the set of permutations in the interval $[u,v]$
that end with $k$.
We show that the following summation
\begin{equation}\label{annj}
\sum_{w\in [u,v]_k}(-1)^{\ell(w)-\ell(u)} \R_{u,w}(q)\R_{w,v}(q),
\end{equation}
 equals zero or a power of $q$
up to a sign. More precisely,  let $D(u,v)$ denote the set of indices where $u$ and $v$ differ,
that is,
\[D(u,v)=\{1\leq i\leq n\colon u(i)\neq v(i)\}.\]
Suppose that $D(u,v)=\{i_1,i_2,\ldots,i_j\}_<$, that is, $D(u,v)=\{i_1, i_2, \ldots, i_j\}$ and $i_1<i_2<\cdots<i_j$.
Let $b_{1}<b_{2}<\cdots<b_{j}$ be the  values of
$u(i_1),u(i_2),\ldots,u(i_j)$ listed in increasing order.
We say that $[u,v]$ is an S-interval if it satisfies the following conditions:
\begin{itemize}
\item[(1)]  $i_j=n$ and $u(i_j)=b_j$;

\item[(2)] the values in
$\{b_1,b_2,\ldots,b_j\}$ that are
greater than $u(i_1)$ appear increasingly in $u$, whereas
the values in
$\{b_1,b_2,\ldots,b_j\}$ that are
less than $u(i_1)$ appear decreasingly in $u$. In other words,
if $b_{j_0}=u(i_1)$ ($1\leq j_0<j$), then the subsequence $u(i_{2})\cdots u(i_{j-1})$ of $u$ is a
shuffle of the decreasing sequence $b_{j_0-1}\cdots b_1$
and the increasing sequence $b_{j_0+1}\cdots b_{j-1}$.

\item[(3)]
$v$ is obtained from $u$ by rotating  the values $b_1,b_2,\ldots,b_j$.
That is, for $m\not\in D(u,v)$ we have $v(m)=u(m)$, and for $m\in D(u,v)$,
if $m=i_j$ then $v(m)=b_1$, and if  $m\neq i_j$ then we have
$v(m)=b_{t+1}$, where $t$ is the index such that $u(m)=b_t$.

\end{itemize}

For example, let $u=4\textbf{3}\textbf{2}\textbf{5}9\textbf{6}\textbf{1}\textbf{7}\textbf{8}$ and
$v=4\textbf{5}\textbf{3}\textbf{7}9\textbf{2}\textbf{6} \textbf{8}\textbf{1}$.
So we have  $D(u,v)=\{2,3,4,6,7,8,9\}$ and $(b_1,\ldots,b_7)=(1,2,3,5,6,7,8)$. The elements $b_i$ are written in boldface.
It is easily checked that $u$ and $v$ satisfy Conditions (1), (2) and (3).
 Notice that the subsequence $u(3)u(4)u(6)u(7)u(8)=25617$ is a shuffle of $21$ and $567$. Thus $[u,v]$ is
an S-interval.

Our main theorem of this section can be stated as follows.

\begin{thm}\label{refine}
Assume that $u<v\in S_n$. Let  $m$ be the smallest index such that $u(m)\neq v(m)$.
If $[u,v]$ is  an S-interval and $k=u(m)$ or $k=v(m)$, then we have
\begin{equation}\label{4-2}
\sum_{w\in [u,v]_k}(-1)^{\ell(w)-\ell(u)}
\widetilde{R}_{u,w}(q) \widetilde{R}_{w,v}(q)=
      (-1)^{r}q^{s-1},
\end{equation}
where $s=|D(u,v)|$ and  $r=|\{j\in D(u,v)\colon \, u(j)>k\}|$; Otherwise, we have
\begin{equation}
\sum_{w\in [u,v]_k}(-1)^{\ell(w)-\ell(u)}
\widetilde{R}_{u,w}(q) \widetilde{R}_{w,v}(q)=0.
\end{equation}
\end{thm}

We give an example before presenting a proof of Theorem \ref{refine}.
Let $n=7$, $u=2354167$ and  $v=3564271$.
It is easy to check that $[u,v]$ is an S-interval.
Consider the sum \eqref{annj} over $[u,v]_2$,  $[u,v]_3$ and $[u,v]_5$, where
\begin{align*}
& [u,v]_2=\{3456172,3465172,3546172,3564172\},\\
& [u,v]_3=\{2456173,2465173,2546173,2564173\},\\
& [u,v]_5=\{2364175,2463175,3264175,3462175\}.
\end{align*}
It can be checked that
\[
\begin{array}{ll}
\widetilde{R}_{2354167,\,3456172}(q)=q^4,\hspace{1.5cm} & \widetilde{R}_{3456172,\,3564271}(q)=q^3,\\[3pt]
\widetilde{R}_{2354167,\,3465172}(q)=q^5, & \widetilde{R}_{3465172,\,3564271}(q)=q^2,\\[3pt]
\widetilde{R}_{2354167,\,3546172}(q)=q^5, & \widetilde{R}_{3546172,\,3564271}(q)=q^2,\\[3pt]
\widetilde{R}_{2354167,\,3564172}(q)=q^6+q^4, & \widetilde{R}_{3564172,\,3564271}(q)=q,
\end{array}\]
and thus we have
\[\sum_{w\in [u,v]_2}(-1)^{\ell(w)-\ell(u)}
\widetilde{R}_{u,w}(q)\widetilde{R}_{w,v}(q)=q^5.\]
Similarly, we have
\[\sum_{w\in [u,v]_3}(-1)^{\ell(w)-\ell(u)}
\widetilde{R}_{u,w}(q)\widetilde{R}_{w,v}(q)=-q^5.\]
Since
\[
\begin{array}{ll}
\widetilde{R}_{2354167,\,2364175}(q)=q^2,\hspace{1.5cm} & \widetilde{R}_{2364175,\,3564271}(q)=q^5+q^3,\\[3pt]
\widetilde{R}_{2354167,\,2463175}(q)=q^3, & \widetilde{R}_{2463175,\,3564271}(q)=q^4,\\[3pt]
\widetilde{R}_{2354167,\,3264175}(q)=q^3, & \widetilde{R}_{3264175,\,3564271}(q)=q^4+q^2,\\[3pt]
\widetilde{R}_{2354167,\,3462175}(q)=q^4, & \widetilde{R}_{3462175,\,3564271}(q)=q^3,
\end{array}\]
we find that
\[\sum_{w\in [u,v]_5}(-1)^{\ell(w)-\ell(u)}
\widetilde{R}_{u,w}(q)\widetilde{R}_{w,v}(q)=0.\]

To prove  Theorem \ref{refine}, let $P_k(u,v)$ denote the set of
 $V$-paths  from
$u$ to $v$  with  bottom $w$
contained in $[u,v]_k$. We shall give an involution $\widetilde{I}$
on $P_k(u,v)$. From now on, assume that $\prec$ is the  lexicographic ordering, namely, the transpositions are ordered by
\[(1,2)\prec(1,3)\prec\cdots\prec(1,n)\prec(2,3)\prec\cdots\prec(n-1,n).\]
We also write   a Bruhat path  $\Delta=w_0\rightarrow w_1\rightarrow\cdots \rightarrow w_r$
as
\begin{equation*}
\Delta=(t_1,\,t_2,\ldots, t_r),
\end{equation*}
where $w_i=w_{i-1}\, t_i$ for $1\leq i\leq r$. It is clear that $w_i=w_0\,t_1\cdots t_i$.

\noindent
\textbf{The Involution $\widetilde{I}$:}
Let $(\Delta_1,\Delta_2)$ be a
$V$-path in $P_k(u,v)$ with bottom $w$.
Write
\[\Delta_1=u_0\rightarrow u_1\rightarrow\cdots \rightarrow u_i\ \ \text{and}\ \
\Delta_2=v_0\rightarrow v_1\rightarrow\cdots \rightarrow v_j,\]
where $u_0=u$, $u_i=v_0=w$ and $v_j=v$.
Assume that $u_r=u_{r-1}\,t_r$ for $1\leq r\leq i$, and
$v_s=v_{s-1}\,t_s'$ for $1\leq s\leq j$. Since
$(\Delta_1,\Delta_2)$ is a $V$-path, we have
\[t_1\succ t_2\succ\cdots\succ t_i\ \ \ \text{and}\ \ \ t_1'\prec t_2'\prec \cdots\prec t_j'.\]

Let $t=\min\{t_i,t_1'\}$.
To define $\widetilde{I}$, here are three cases.

\noindent
Case 1. $t$ is an internal transposition, that is, $t=(a,b)$ and $1\leq a<b<n$.
Set
\[\widetilde{I}(\Delta_1,\Delta_2)=I(\Delta_1,\Delta_2).\]

In this case,
the bottom node of $\widetilde{I}(\Delta_1,\Delta_2)$ is also contained
 in $[u,v]_k$. Thus  $\widetilde{I}(\Delta_1,\Delta_2)$ remains  a $V$-path
 in $P_k(u,v)$.

\noindent
Case 2. $t$ is a boundary transposition, that is, $t=(a,n)$ for some $a<n$, and  there is at least one  internal transposition
among the transpositions
$t_1,\ldots, t_i, t_1',\ldots, t_j'$.
Let $\widetilde{t}$ be the smallest internal transposition among
 $t_1,\ldots, t_i, t_1',\ldots, t_j'$.
We consider the following two subcases.

Subcase 1. $\widetilde{t}$
belongs to $\{t_1,\ldots, t_i\}$.
Assume that $t_{i_0}=\widetilde{t}$, where $1\leq i_0\leq i$.
If the transpositions $t_1',\ldots,t_j'$ are all smaller than $\widetilde{t}$,
then set
\begin{equation*}
\Delta_1'=(t_1,\,\ldots,\,\widehat{t_{i_0}},\,\ldots,\, t_i)\ \ \text{and}\ \
\Delta_2'=(t_1',\,t_2',\ldots, t_j',\,t_{i_0}),
\end{equation*}
where $\widehat{t_{i_0}}$ indicates that the element $t_{i_0}$ is missing.
Define $\widetilde{I}(\Delta_1,\Delta_2)=(\Delta_1',\Delta_2')$.

We claim that $(\Delta_1',\Delta_2')$ is contained in $P_k(u,v)$.
By the choice of $t_{i_0}$,
we see that the boundary transpositions $t_{i_0+1},\ldots,t_i$   commute with $t_{i_0}$.
Thus the bottom of  $(\Delta_1',\Delta_2')$ is given by
\[u\,t_1\cdots t_{i_0-1}\widehat{t_{i_0}}t_{i_0-1}\cdots t_i
= u\,t_1\cdots t_{i_0-1}t_{i_0}t_{i_0+1}\cdots t_i\,t_{i_0}=w\,t_{i_0}.\]
Clearly, $w\,t_{i_0}$ is a permutation in $[u,v]_k$. Hence the claim is proved.

If there exists a transposition among $t_1',\ldots,t_j'$ that is greater than $\widetilde{t}$, let $1\leq p\leq j$ be the smallest index such that $t_{p}'\succ\widetilde{t}$.
Set
\[
\Delta_1'=(t_1,\,\ldots,\,\widehat{t_{i_0}},\,\ldots,\, t_i)\ \ \text{and}\ \
\Delta_2'=(t_1',\,\ldots,\,t_{p-1}',\,t_{i_0},\,t_{p}',\,\ldots,\, t_j').
\]
Define
$\widetilde{I}(\Delta_1,\Delta_2)=(\Delta_1',\Delta_2')$.

Similarly, it can be checked that
the bottom of $(\Delta_1',\Delta_2')$ is $w\,t_{i_0}$. Hence $(\Delta_1',\Delta_2')$ is a $V$-path
in $P_k(u,v)$.

Subcase 2. $\widetilde{t}$
belongs to $\{t_1',\ldots, t_j'\}$.
Assume that $t_{j_0}'=\widetilde{t}$, where $1\leq j_0\leq j$.

If $t_1,\ldots,t_i$ are all smaller than $\widetilde{t}$,
then set
\[
\Delta_1'=(t_{j_0}',\,t_1,\,\ldots, t_i)\ \ \text{and}\ \
\Delta_2'=(t_1',\,\ldots,\,\widehat{t_{j_0}'},\,\ldots,\, t_j').
\]
Define
$\widetilde{I}(\Delta_1,\Delta_2)=(\Delta_1',\Delta_2')$.

If there exists a transposition among $t_1,\ldots,t_i$ that is greater than $\widetilde{t}$, let $1\leq p\leq i$ be the largest  index such that $t_{p}\succ \widetilde{t}$.
Set
\[
\Delta_1'=(t_1,\,\ldots,\,t_{p},\,t_{j_0}',\,t_{p+1},\,\ldots,\, t_i)\ \ \text{and}\ \
\Delta_2'=(t_1',\,\ldots,\,\widehat{t_{j_0}'},\,\ldots,\, t_j').
\]
Define $\widetilde{I}(\Delta_1,\Delta_2)=(\Delta_1',\Delta_2')$.
Using the argument in subcase 1, it can be seen  that
the bottom of $(\Delta_1',\Delta_2')$ is $w\,t_{i_0}$.
Thus $(\Delta_1',\Delta_2')$ is a $V$-path
in $P_k(u,v)$.

\noindent
Case 3. The
transpositions $t_1,\ldots, t_i, t_1',\ldots, t_j'$
are all boundary transpositions.
In this case, set $(\Delta_1,\Delta_2)$ to be
a fixed point of $\widetilde{I}$, that is,
$\widetilde{I}(\Delta_1,\Delta_2)=(\Delta_1,\Delta_2)$.

The construction of $\widetilde{I}$ implies the following property.

\begin{thm}\label{r-9}
The map $\widetilde{I}$ is a length preserving involution on    $P_k(u,v)$.
Moreover, $\widetilde{I}$ reverses the sign of $(\Delta_1,\Delta_2)$
unless $(\Delta_1,\Delta_2)$ is a fixed point.
\end{thm}

To finish the proof of Theorem \ref{refine}, we  need the following lemma.

\begin{lem}\label{lemma-I}
The involution $\widetilde{I}$ has at most one fixed point in $P_k(u,v)$.
Moreover, $\widetilde{I}$ has a fixed point in $P_k(u,v)$
if and only if
$[u,v]$ is an $S$-interval and $k=u(m)$ or $k=v(m)$, where $m$ is the
smallest index such that $u(m)\neq v(m)$.
\end{lem}

\pf We first show that $\widetilde{I}$ has at most one fixed point in $P_k(u,v)$.
Suppose that $(\Delta_1,\Delta_2)\in P_k(u,v)$ is a $V$-path   with bottom $w$ that is fixed by $\widetilde{I}$.
Let
\begin{equation*}
\Delta_1=(t_1,\,t_2,\ldots, t_i)\ \ \text{and}\ \ \Delta_2=(t_1',\,t_2',\ldots, t_j').
\end{equation*}
By the construction of $\widetilde{I}$, we see that
$t_1,\ldots, t_i$ and $t_1',\ldots, t_j'$ are all boundary transpositions.
Assume that
\begin{equation*}
t_1=(p_1,n),\,\ldots\, ,t_i=(p_i,n)\ \text{  and}\ \ t_1'=(p_1',n),\,\ldots\, , t_j'=(p_j',n).
\end{equation*}
Since $\Delta_1$ and $\Delta_2$ are Bruhat paths, we see that
\begin{equation}\label{r-3}
u(n)>u(p_1)>\cdots>u(p_i)=k=w(n)>w(p_1')>\cdots>w(p_j').
\end{equation}
Moreover, noting that $t_1\succ t_2\succ\cdots\succ t_i$ and $t_1'\prec t_2'\prec \cdots\prec t_j'$,
we have
\begin{equation*}
n>p_1>\cdots>p_i\ \ \text{ and}\ \ p_1'<\cdots<p_j'<n.
\end{equation*}

By \eqref{r-3} and relation $w=u\,(p_1,n)\cdots (p_i,n)$, it is easily seen that
$\{p_1,\ldots,p_i\}\cap\{p_1',\ldots,p_j'\}=\emptyset$. This implies that
\begin{equation*}
w(p_1')=u(p_1'),\,\ldots\, , w(p_j')=u(p_j'),
\end{equation*}
and so \eqref{r-3} becomes
\begin{equation}\label{r-4}
u(n)>u(p_1)>\cdots>u(p_i)=k=w(n)>u(p_1')>\cdots>u(p_j').
\end{equation}

Observing that
$
\{p_1,\ldots,p_i\}\cup\{p_1',\ldots,p_j'\}\cup\{n\}=D(u,v)$,
in view of \eqref{r-4} we see that once $u,\,v$ and $k$ are given,
the values $i,j$ as well as $p_1,\ldots, p_i, p_1',\ldots, p_j'$ are uniquely determined.
This leads to the conclusion that $\widetilde{I}$ has at most one fixed point in $P_k(u,v)$.

It remains to prove the assertion that $P_k(u,v)$ has a fixed point under $\widetilde{I}$
if and only if
$[u,v]$ is an S-interval and $k=u(m)$ or $k=v(m)$.
By the above argument, we find that if $\widetilde{I}$ has a fixed point  in $P_k(u,v)$, then $[u,v]$ is an S-interval and $k=u(p_i)=v(p_1')$.
Since $m=\min\{p_i,p_1'\}$, we obtain that $k=u(m)$ if $p_i<p_1'$ and $k=v(m)$ if $p_i>p_1'$.
Conversely, if $[u,v]$ is an S-interval, it is easy to construct a
$V$-path in $P_k(u,v)$ fixed by $\widetilde{I}$, where $k=u(m)$ or $k=v(m)$.
This completes the proof.
\qed

We give an example to illustrate the above lemma.
Let $n=9$, $u=432596178$ and $v=453792681$.
It can be seen that $[u,v]$ is an S-interval.
The fixed point in $P_3(u,v)$ is $(\Delta_1,\Delta_2)$
with $\Delta_1=((8,9),(6,9),(4,9),(2,9))$ and $\Delta_2=((3,9),(7,9))$.
The fixed point in $P_5(u,v)$ is $(\Delta_3,\Delta_4)$
with $\Delta_3=((8,9),(6,9),(4,9))$ and $\Delta_4=((2,9)(3,9),(7,9))$.

We are now ready to complete the proof of Theorem
\ref{refine}.

\noindent
\textit{Proof of Theorem \ref{refine}.}
By Theorem \ref{r-9} and Lemma \ref{lemma-I},
we only need to consider the case when $[u,v]$ is an $S$-interval and $k=u(m)$ or
$k=v(m)$. In this case, we have
\[\sum_{w\in [u,v]_k}(-1)^{\ell(w)-\ell(u)}
\widetilde{R}_{u,w}(q)\widetilde{R}_{w,v}(q)=
(-1)^{\ell(\Delta_1)}q^{\ell(\Delta_1)+\ell(\Delta_2)},\]
where $(\Delta_1,\Delta_2)$ is the unique $V$-path
 in $P_k(u,v)$ that is fixed by $\widetilde{I}$.
Evidently, we have
\[\ell(\Delta_1)+\ell(\Delta_2)=|D(u,v)|-1.\]
Moreover, it is easy to check that
\[\ell(\Delta_1)=|\{j\in D(u,v)\colon \, u_j>k\}.\]
This completes the proof.
\qed

\vspace{.2cm} \noindent{\bf Acknowledgments.}
This work was
supported by the 973 Project, the PCSIRT Project of the Ministry of
Education, and the National Science Foundation of China.

\end{document}